\documentclass[a4paper,11pt]{article}
\usepackage{amsmath}
\usepackage{amssymb}
\usepackage{enumerate}
\usepackage{theorem}
\usepackage{array}
\usepackage{bm}
\newtheorem{theorem}{Theorem}[section]
\newtheorem{lemma}[theorem]{Lemma}
\newtheorem{corollary}[theorem]{Corollary}
\newtheorem{proposition}[theorem]{Proposition}
\theorembodyfont{\rmfamily}

\newtheorem{remark}[theorem]{Remark}
\newtheorem{example}[theorem]{Example}

\newtheorem{notationinintro}{Notation}

\newcommand{\proof}{\noindent \mbox{\em Proof.\hspace*{2mm}}}
\newcommand{\qed}{\hfill \mbox{$  \Box $}}

\DeclareMathOperator{\rank}{rank}

\newcommand{\gyokan}{\vskip 6pt}
\title{Fragments of plane filling curves of degree $q+2$ over the finite field of $q$ elements,\\
and of affine-plane filling curves of degree $q+1$}
\author{
Masaaki Homma
\\
 Department of Mathematics and Physics\\
Kanagawa University\\
Hiratsuka 259-1293, Japan\\
homma@kanagawa-u.ac.jp
}

\date{}

\begin{document}
\maketitle
\begin{abstract}
Nonsingular plane curves over a finite field $\mathbb{F}_q$ of degree $q+2$
passing through all the $\mathbb{F}_q$-points of the plane admit
a representation by $3\times 3$ matrices over $\mathbb{F}_q$.
We classify their degenerations by means of the matrix representation.
We also discuss the similar problem for the affine-plane filling curves of degree $q+1$.
\end{abstract}

\section{Introduction}
This paper is a continuation of \cite{hom-kim2013} published some years ago.
As is used in the title of both papers, a plane filling curve is a plane curve over $\mathbb{F}_q$ passing through all $\mathbb{F}_q$-points of the projective plane $\mathbb{P}^2$.
The set of 
all $\mathbb{F}_q$-points of the plane will be denoted by $\mathbb{P}^2(\mathbb{F}_q)$.
In the previous paper, we showed that the minimum degree of nonsingular plane filling curves over $\mathbb{F}_q$ is $q+2$, and
described those curves of degree $q+2$.
More precisely, to each plane filling curve of degree $q+2$,
there is a $3\times 3$-matrix $A$ over $\mathbb{F}_q$ so that
the plane curve is defined by
\begin{equation}\label{equationA}
(x,y,z)A \,{}^t\! (y^q z - y z^q, z^q x - z x^q, x^q y - x y^q)=0.
\end{equation}
(Here, we make a correction to \cite[p.970, Introduction]{hom-kim2013};
in the expression of the degree $q+2$ homogeneous part of the ideal of
$\mathbb{P}^2(\mathbb{F}_q)$, ``$A \in GL(3, \mathbb{F}_q)$" should be read as
``$A \in \mathbb{M}(3, \mathbb{F}_q)$",
where $\mathbb{M}(3, \mathbb{F}_q)$ is the set of $3\times 3$ matrices
over $\mathbb{F}_q$, however, the expression in \cite{hom-kim2013} is also true
if $q>3$ because
$(x,y,z)\mu E \,{}^t\! (y^q z - y z^q, z^q x - z x^q, x^q y - x y^q)$
is the zero polynomial
for $\mu \in \mathbb{F}_q$ and the identity matrix $E$.)
We gave a condition on such a plane filling curve to be nonsingular
in terms of the corresponding matrix.
Happily any irreducible plane filling curve of degree $q+2$ is nonsingular,
and irreducible plane filling curves of degree $q+2$ had been studied well
by Tallini \cite{tal2} more than half a century ago. 
In Section~2, we will review some parts of \cite{hom-kim2013}
including these matters plainly.

The purpose of this paper is to study
reducible plane filling curves of degree $q+2$.
Precisely we want to know how they split into irreducible pieces,
which will be done in Sections~3 and 4.
A remarkable phenomenon is that fragments other than $\mathbb{F}_q$-lines of the splits are only four kinds of curves,
each of which is a particular curve in the sense of the Sziklai bound;
three of the four are maximal curves
of degrees $q+1$, $q$ and $q-1$ respectively,
and the remaining one is an affine-plane filling curve, which attains the second maximum for degree $q+1$ curves.

The main theorem of the previous paper was that
\textit{the curve defined by $(\ref{equationA})$ is nonsingular if and only if
the characteristic polynomial of $A$ is irreducible over $\mathbb{F}_q$.}
On the same lines, we will show that
\textit{the curve defined by $(\ref{equationA})$ has a non-linear component
if and only if
the characteristic polynomial of $A$ is also the minimal polynomial.}

In the last two sections, we will study the ``degeneration" of the family of affine-plane filling curves of degree $q+1$.
Let us fix an $\mathbb{F}_q$-line, say the line $l_{\infty}$ defined by
$z=0$, and let
$\mathbb{A}^2=\mathbb{P}^2 \setminus l_{\infty}$.
A plane curve $C$ over $\mathbb{F}_q$ is an affine-plane filling curve
if $C(\mathbb{F}_q)= \mathbb{A}^2(\mathbb{F}_q)$.
The minimum degree of affine-plane curves is $q+1$,
and such a curve is represented by a $2\times 3$ matrix over $\mathbb{F}_q$
with a certain condition (see, Theorem~\ref{affineplanefilling} below).
Without the condition, the $2\times 3$ matrix gives a possibly reducible curve which still satisfies $C(\mathbb{F}_q) \supset \mathbb{A}^2(\mathbb{F}_q)$.
The problem is, again, how such curves split into irreducible pieces.
We will completely describe the splits by means of those $2\times 3$ matrices.

\begin{notationinintro}
Throughout this paper, we fix a finite field $\mathbb{F}_q$
of $q$-elements.
The multiplicative group of  $\mathbb{F}_q$
is denoted by
$\mathbb{F}_q^{\ast}$. 
For two positive integers $m$ and $n$,
$\mathbb{M}(m\times n, \mathbb{F}_q)$
denotes the set of $m\times n$ matrices.
When $m=n$, $\mathbb{M}(m\times m, \mathbb{F}_q)$
is simply denoted by $\mathbb{M}(m, \mathbb{F}_q)$.

We also fix a projective plane $\mathbb{P}^2$
with homogeneous coordinates $x, y, z$
defined over $\mathbb{F}_q$.
When $F =F(x,y,z) \in \mathbb{F}_q[x,y,z]$ is homogeneous,
the symbol $\{F=0\}$ denotes
the algebraic curve in $\mathbb{P}^2$ defined by the equation
$F=0$.
Additionally, $C(\mathbb{F}_q)$ denotes the set of $\mathbb{F}_q$-
points of a curve $C$ over $\mathbb{F}_q$,
and $N_q(C)$ denotes the cardinality of $C(\mathbb{F}_q)$.
\end{notationinintro}
\section{Plane filling curve of degree $q+2$}

It is well known that the ideal $\mathfrak{h}$
of the set of all $\mathbb{F}_q$-points $\mathbb{P}^2(\mathbb{F}_q)$
in $\mathbb{F}_q[x,y,z]$ is generated by
\[
  \begin{array}{ccl}
  U &:=& y^q z - y z^q, \\
  V &:=& z^q x - z x^q, \\
  W &:=& x^q y - x y^q .
  \end{array}
\]
For convenience, we also use the notation
$U(y,z)$, $V(z,x)$ and $W(x,y)$ for $U, V$ and $W$ respectively.
In particular, the degree of a (reducible) plane filling curve is at least $q+1$, and any one of degree $q+1$
is a union of $q+1$ $\mathbb{F}_q$-lines passing through a certain point
\cite[Proposition 2.3]{hom-kim2013}.

The next degree, namely $q+2$, is of great interest, which is the minimum degree of existence of (absolutely) irreducible plane filling curves.
The vector space of homogeneous elements of degree $q+2$ of $\mathfrak{h}$
is denoted by $\mathfrak{h}_{q+2}$.
Obviously, any element of $\mathfrak{h}_{q+2}$
can be represented as
\[
F_A(x,y,z) := (x,y,z) A 
\begin{pmatrix} U\\ V \\ W\end{pmatrix} 
\]
for a certain matrix $A \in \mathbb{M}( 3, \mathbb{F}_q).$
\begin{lemma}\label{kernel}
The kernel of the surjective map
\[
\mathbb{M}( 3, \mathbb{F}_q) \ni A 
       \mapsto F_A(x,y,z) \in \mathfrak{h}_{q+2}
\]
is $\{\mu E \mid \mu \in \mathbb{F}_q \}.$
\end{lemma}
\proof
See \cite[Lemmma 3.1]{hom-kim2013}.
\qed

\gyokan

For $A \in \mathbb{M}( 3, \mathbb{F}_q)
     \setminus \{\mu E \mid \mu \in \mathbb{F}_q \},$
$C_A$ denotes the curve
defined by $ F_A(x,y,z)=0,$
which is a plane filling curve of degree $q+2$.

\begin{lemma}\label{howtotransformB}
Let us consider the coordinate change by
an $\mathbb{F}_q$-linear transformation
\begin{equation}\label{lineartransformB}
\begin{pmatrix}
x \\ y \\ z
\end{pmatrix}
= B
\begin{pmatrix}
x' \\ y' \\ z'
\end{pmatrix}
\text{with \ }
B = (b_{ij}) \in GL(3, \mathbb{F}_q).
\end{equation}
Let
$U' =U(y',z')$, $V'=(z',x')$ and $W'=W(x',y')$.
Then the relation
\[
\begin{pmatrix}
U(b_{21}x' + b_{22}y' + b_{23}z', b_{31}x' + b_{32}y' + b_{33}z')\\
V(b_{31}x' + b_{32}y' + b_{33}z', b_{11}x' + b_{12}y' + b_{13}z')\\
W(b_{11}x' + b_{12}y' + b_{13}z', b_{21}x' + b_{22}y' + b_{23}z')
\end{pmatrix}
= (\det B) {}^t\! B^{-1}
\begin{pmatrix}
U'\\ V' \\ W'
\end{pmatrix}
\]
holds.
\end{lemma}
\proof
See \cite[Lemmma 2.2]{hom-kim2013}.
\qed

\gyokan

This lemma implies that the polynomial
$F_A(x,y,z)$ in $x,y,z$ is transformed into 
$F_{(\det B){}^t\!BA{}^t\!B^{-1}}(x',y',z')$ in $x',y',z'$
by the linear transformation (\ref{lineartransformB}).

\begin{corollary}\label{corollaryequivalent}
Let $A$ and $A'$ in $\mathbb{M}( 3, \mathbb{F}_q) \setminus \{\mu E \mid \mu \in \mathbb{F}_q \}.$
Two plane curves $C_A$ and $C_{A'}$ are projectively equivalent over $\mathbb{F}_q$ if and only if
there are $B \in GL(3, \mathbb{F}_q)$, $\rho \in \mathbb{F}_q^{\ast}$
and $\mu \in \mathbb{F}_q$ such that
\begin{equation}\label{projectivelyequiv}
A' = \rho {}^t\!BA{}^t\!B^{-1} + \mu E.
\end{equation}
\end{corollary}
\proof
Suppose $C_A$ and $C_{A'}$ are projectively equivalent to each other over $\mathbb{F}_q$.
Then there is a projective transformation $B \in GL(3, \mathbb{F}_q)$
such that
$F_{A'} =\rho' F_{(\det B){}^t\!BA{}^t\!B^{-1}}$
for some $\rho' \in \mathbb{F}_q^{\ast}$
by Lemma~\ref{howtotransformB}.
Hence $A' - \rho'(\det B){}^t\!BA{}^t\!B^{-1}= \mu E$
for some $\mu \in \mathbb{F}_q$
by Lemma~\ref{kernel}.
Conversely, suppose $A' = \rho {}^t\!BA{}^t\!B^{-1} + \mu E$.
Then
$F_{A'} =  0$ and $F_{(\det B){}^t\!BA{}^t\!B^{-1}}=0$  define
the same curve $C_{A'}$ by Lemma~\ref{kernel},
and hence $C_A$ and $C_{A'}$ are projectively equivalent  over $\mathbb{F}_q$
by Lemma~\ref{howtotransformB}.
\qed

\gyokan

For a matrix $A \in \mathbb{M}(3 , \mathbb{F}_q)$,
$f_A(t)$ denotes the characteristic polynomial $|tE-A|$ of $A$.
Note that if there is the relation~(\ref{projectivelyequiv})
between $A$ and $A'$, then
$f_{A'}(t) = \rho^3 f_A\left( \frac{t-\mu}{\rho}\right).$

\gyokan

We supplemented Tallini's work \cite{tal1, tal2}
by proving any plane filling irreducible curve of degree $q+2$
to be nonsingular.
The following theorem is essentially a reformation of
what we have shown in \cite{hom-kim2013}.

\begin{theorem}[\cite{tal2} and \cite{hom-kim2013}]\label{previousmaintheorem}
For $A \in \mathbb{M}(3 , \mathbb{F}_q)$,
the following conditions are equivalent{\rm :}
\begin{enumerate}[{\rm (a)}]
\item $C_A$ is nonsingular{\rm ;}
\item $C_A$ is absolutely irreducible{\rm;}
\item $C_A$ is irreducible over $\mathbb{F}_q${\rm;}
\item $C_A$ has no $\mathbb{F}_q$-linear components{\rm ;}
\item $C_A$ is nonsingular at each $\mathbb{F}_q$-point{\rm ;}
\item the characteristic polynomial $f_A(t)$ is irreducible
over $\mathbb{F}_q$.
\end{enumerate}
\end{theorem}
\proof
The equivalence (a)$\Leftrightarrow$(f)  was showen
at \cite[Theorem 3.2]{hom-kim2013}.
Other equivalences from (b) to (f) but (d) are due to \cite{tal2}.
Actually, although he didn't use the matrix representation for those curves,
the condition on coefficients of his canonical forms
of irreducible filling curves agrees with the above (f).

Another approach to show the equivalence (a)$ \Leftrightarrow$(b) 
is to use a property of the automorphism group
of an irreducible filling curve, namely, it has the Singer automorphism,
which was also mentioned in \cite[the end of Section~5]{hom-kim2013}
(see, also \cite[Theorem 2.1]{dur2018}).

The implication (d)$\Rightarrow$(c) comes from the Sziklai bound
(\ref{sziklai}), which will be reviewed in the next section.
Actually, suppose $C_A$, which has $q^2+q+1$ $\mathbb{F}_q$-points, decomposes into
two curves $C_A = C_1 \cup C_2$ over $\mathbb{F}_q$,
but each of them has no $\mathbb{F}_q$-linear components.
Unless $q = 4$ and one of the $C_i$'s is the exceptional curve (\ref{exception}) in the next section,
$N_q(C_i) \leq (\deg C_i -1)q +1$.
Since
$\deg C_1 + \deg C_2 = q+2$,
\begin{align*}
N_q(C_A) &\leq (\deg C_1 + \deg C_2 -2)q +2\\
         &= q^2 +2 < q^2 + q + 1,
\end{align*}
which is absurd.
When $q=4$ and a component of $C_A$ is the exceptional curve (\ref{exception}),
which has $14$ $\mathbb{F}_4$-points, then other one is an irreducible conic,
which has $5$ $\mathbb{F}_4$-points. Hence
$N_4(C_A) \leq 19 < N_4(\mathbb{P}^2)= 21,$
which is also absurd.
\qed

\gyokan

\begin{remark}
Recently, Duran Cunha \cite[Theorem 2.2]{dur2018}
proved that irreducible plane filling curves of degree $q+2$
over $\mathbb{F}_q$ are projectively equivalent to each other
over the algebraic closure of $\mathbb{F}_q$,
precisely, those curves are projectively equivalent to
the curve defined by
\[
xy^{q+1} + yz^{q+1} + zx^{q+1} =0
\]
over $\mathbb{F}_{q^{3(q^2 + q +1)}}.$
In his proof, he used only the condition~(f)
in Theorem~\ref{previousmaintheorem}.
Actually the cubic equation (2) in his paper \cite{dur2018}
is exactly the characteristic polynomial of
the plane filling curve (1) in \cite{dur2018}
in our context.
So this also gives the third approach to show that any irreducible plane filling curve of degree $q+2$ is nonsingular.

\end{remark}
\section{Plane curves with many points of non-small degrees}
Reducible plane filling curves of degree $q+2$
must contain an $\mathbb{F}_q$-line as an irreducible component
because of Theorem~\ref{previousmaintheorem} (d)$\Rightarrow$(c).
We want to know what curves other than lines appear as components
of reducible plane filling curves.
Those curves still have many $\mathbb{F}_q$-points with respect to their degrees.
We list those curves in advance.

The Sziklai bound, which was already used in the previous section, says that
if a degree $d$ curve $C$ has no $\mathbb{F}_q$-linear components,
then
\begin{equation}\label{sziklai}
N_q(C) \leq (d-1)q+1
\end{equation}
except for the case where $d=q=4$ and $C$ is projectively equivalent to
the curve defined by the equation
\begin{equation}\label{exception}
(x+y+z)^4+(xy+yz+zx)^2+xyz(x+y+z) =0.
\end{equation}
For the Sziklai bound, consult \cite{hom-kim2010b}.

\begin{example}\label{degreeqplusone}
When $d=q+1$, the curve
\begin{equation}\label{qplus1}
x^{q+1} - x^2z^{q-1} + y^qz - yz^{q} = 0
\end{equation}
attains the upper bound in (\ref{sziklai}),
and conversely if equality holds in (\ref{sziklai})
for a curve without $\mathbb{F}_q$-linear components,
then the curve is projectively equivalent to (\ref{qplus1})
over $\mathbb{F}_q$
except for the case $q=3$ and $4$ \cite{hom-kim2011}.
The curve (\ref{qplus1}) is obviously nonsingular.
\end{example}
\begin{example}\label{affinefilling}
If $C(\mathbb{F}_q) =
\mathbb{P}^2(\mathbb{F}_q) \setminus (\mbox{an $\mathbb{F}_q$-line}),$
then $\deg C \geq q+1$, and such a curve of degree $q+1$ exists \cite[Theorem 3.5]{hom-kim2018}.
The curve mentioned above is a second maximum curve of degree $q+1$
in the sense of the number of $\mathbb{F}_q$-points.
Explicitly such a curve is projectively equivalent to the curve
defined by the equation
\[
(x^q-xz^{q-1}, y^q-yz^{q-1})
\left(
 \begin{array}{ccc}
   a_0 & a_1 & a_2 \\
   b_0 & b_1 & b_2
 \end{array}
\right)
\left(
\begin{array}{c}
x \\ y \\ z
\end{array}
\right) = 0,
\]
where $a_0, \dots , b_2 $ are elements of $\mathbb{F}_q$
and the polynomial
$
(s,t)
\left(
 \begin{array}{cc}
   a_0 & a_1 \\
   b_0 & b_1 
 \end{array}
\right)
\left(
\begin{array}{c}
s \\ t
\end{array}
\right)
$
in $s$ and $t$ is irreducible over $\mathbb{F}_q$.
Note that this curve is absolutely irreducible, and
has a unique singular point.

The last two sections of this paper will be devoted to investigating the ``degeneration" of those curves in the sense of their coefficient matrices.
\end{example}
\begin{example}\label{degreeq}
Let $C$ be a plane curve over $\mathbb{F}_q$, of degree $q$ and 
without $\mathbb{F}_q$-linear components.
$N_q(C)$ is equal to the upper bound of (\ref{sziklai})
if and only if $C$ is projectively equivalent to the curve
\begin{equation}
x^q - xz^{q-1} + x^{q-1}y - y^q =0,
\end{equation}
which is a nonsingular curve.
This is, of course, the maximum value other than $q=4$,
and when $q=4$, this is the second maximum value.
For this example, consult \cite{hom-kim2012}.
\end{example}
\begin{example}\label{degreeqminusone}
There also exists a nonsingular curve $C$ of degree $q-1$
so that $N_q(C)$ attains the upper bound of (\ref{sziklai}).
Actually
\begin{equation}
\alpha x^{q-1} + \beta y^{q-1} + \gamma z^{q-1}=0
\end{equation}
with $\alpha\beta\gamma \neq 0$ and $\alpha + \beta + \gamma = 0$
\end{example} 
is such a curve \cite{szi2008}.
\section{Fragments of plane filling curves}
As we already explained in Section~2, any plane filling curve is represented
by a matrix $A \in \mathbb{M}( 3, \mathbb{F}_q)$
so that the equation of the curve $C_A$ is 
\[
F_A(x,y,z) := (x,y,z) A 
\begin{pmatrix} U\\ V \\ W\end{pmatrix},
\]
and if the matrix $A$ is replaced by ${}^t\!BA{}^t\!B^{-1},$
then $C_{{}^t\!BA{}^t\!B^{-1}}$ is the image of the projective transfrmation
$
\begin{pmatrix} x\\ y \\ z\end{pmatrix}
\mapsto
B^{-1}
\begin{pmatrix} x\\ y \\ z\end{pmatrix}.
$
The equivalence classes of $\mathbb{M}( 3, \mathbb{F}_q)$
by the action $A\mapsto {}^t\!BA{}^t\!B^{-1}$ are determined by
invariant polynomias of $tE-A$.
(For this matter, consult \cite{gan}, for example.)
Since the size of our matrices is only $3$,
the invariant polynomials are determined by the characteristic polynomial
and the minimal polynomial of $A$.

From Theorem~\ref{previousmaintheorem},
$C_A$ splits into two or more components if and only if
the characteristic polynomial $f_A(t)$ of $A$ is reducible over $\mathbb{F}_q$.

When matrices $A$ and $A' \in \mathbb{M}( 3, \mathbb{F}_q)$
are connected by the relation (\ref{projectivelyequiv}),
we indicate this situation by $A \sim A'$.
As we explained in Corollary~\ref{corollaryequivalent}
and Lemma~\ref{kernel},
if $A \not\sim E$,
then $C_A$ and $C_{A'}$ are projectively equivalent over $\mathbb{F}_q$,
and if $A\sim E$,
then $F_A$ is the zero polynomial.

We study $C_A$ according to the following cases:
\begin{enumerate}[{\em Case} 1]
\item $f_A(t) = (t-\alpha)g(t)$, where
 $\alpha \in \mathbb{F}_q$ and $g(t) = t^2 -(bt+a) \in \mathbb{F}_q[t]$
 is irreducible.
 In this case, the elementary divisors of $tE-A$ are
 $\{ g(t), t-\alpha \}$
\item $f_A(t)=(t-\alpha_1)(t-\alpha_2)(t-\alpha_3)$,
 where three elements $\alpha_1, \alpha_2, \alpha_3$ of $\mathbb{F}_q$
 are distinct each other.
 In this case, the elementary divisors of $tE-A$ are
 $\{ t-\alpha_1,\, t-\alpha_2,\, t-\alpha_3 \}$.
\item $f_A(t)=(t-\alpha)^2(t-\beta)$, where
 $\alpha, \beta \in \mathbb{F}_q$ and $\alpha \neq \beta$
 \begin{enumerate}[{\em Case} {3.}1]
  \item The minimal polynomial of $A$ is $f_A(t)$ itself.
  In this case, the elementary divisors of $tE-A$
  are $\{ (t-\alpha)^2 , t-\beta \}.$
  \item The minimal polynomial of $A$ is
  $(t-\alpha)(t-\beta)$.
  In this case, the elementary divisors of $tE-A$
  are $\{ t-\alpha,\, t-\alpha,\, t-\beta \}.$
 \end{enumerate}
\item $f_A(t) = (t-\alpha)^3$.
 \begin{enumerate}[{\em Case} {4.}1]
  \item The minimal polynomial of $A$ is $f_A(t)$ itself.
  In this case, the elementary divisor of $tE-A$ is
  $(t-\alpha)^3$.
  \item The minimal polynomial of $A$ is $(t-\alpha)^2$.
  In this case, the elementary divisors of $tE-A$ are
  $\{ (t-\alpha)^2, t-\alpha\}$.
  \item The minimal polynomial of $A$ is $t-\alpha$.
  In this case, the elementary divisors of $tE-A$ are
  $\{ t-\alpha,\, t-\alpha,\, t-\alpha \}$.
 \end{enumerate}
\end{enumerate}
\begin{theorem}
 Suppose that a plane filling curve $C_A$ is reducible.
 \begin{enumerate}[{\rm (i)}]
 \item In Case1, $C_A$ splits into an $\mathbb{F}_q$-line and an irreducible curve $C_A'$ of degree $q+1$.
 The irreducible curve $C_A'$ is one described in Example~\ref{affinefilling}.
 \item In Case 2, $C_A$ splits into three $\mathbb{F}_q$-lines and
 irreducible curve $C_A'$ of degree $q-1$.
 Those three lines are not concurrent, and the irreducible curve $C_A'$ is one
 described in Example~\ref{degreeqminusone}.
 \item In Case 3.1, $C_A$ splits into two $\mathbb{F}_q$-lines and
 irreducible curve $C_A'$ of degree $q$.
 The irreducible curve $C_A'$ is one described in Example~\ref{degreeq}.
 \item In Case 3.2, $C_A$ splits into $q+2$ $\mathbb{F}_q$-lines.
 Among them, $q+1$ lines pass through a common point, and the remaining one
 does not.
 \item In Case 4.1, $C_A$ splits into an $\mathbb{F}_q$-line and an irreducible curve $C_A'$ of degree $q+1$.
 The irreducible curve $C_A'$ is described in Example~\ref{degreeqplusone}.
 \item In Case 4.2, $C_A$ splits into $q$ reduced $\mathbb{F}_q$-lines
 and one double line over $\mathbb{F}_q$, each of which {\rm (}including the double line{\rm )}, passes through a common point.
 \item In Case 4.3, $F_A$ is the zero polynomial.
 \end{enumerate}
\end{theorem}
\proof
In Case~1, since
$A \sim 
\begin{pmatrix}
0&a&0 \\
1&b&0 \\
0&0&\alpha
\end{pmatrix},$
we may assume that
\begin{align*}
F_A(x,y,z) &= y(y^qz-yz^q) + (ax+by)(z^qx-zx^q) + \alpha z (x^qy -xy^q)\\
    &= z G(x,y,z),
\end{align*}
where
\begin{align*}
G(x,y,z)&= (y^{q+1} -y^2z^{q-1}) +(ax+by)(z^{q-1}x-x^q)
                  +  \alpha  (x^qy -xy^q)\\
        &=(x^q - xz^{q-1}, y^q - yz^{q-1})
           \begin{pmatrix}
            -a & \alpha -b & 0 \\
            -\alpha& 1 &0
           \end{pmatrix}
           \begin{pmatrix}
            x \\ y\\ z
           \end{pmatrix}.
\end{align*}
Since
\[
(s, t)\begin{pmatrix}
            -a & \alpha -b  \\
            -\alpha& 1 
           \end{pmatrix}
               \begin{pmatrix}
            s \\ t
           \end{pmatrix}
  = -as^2 -bst +t^2 = s^2 g(t/s),
\]
it is irreducible by the assumption on the form of $f_A(t)$.

In the case~2,
since
$A \sim 
\begin{pmatrix}
\alpha_1& & \\
 &\alpha_2& \\
 & &\alpha_3
\end{pmatrix},$
we may assume that
\begin{align*}
F_A(x,y,z) &= \alpha_1 x (y^qz-yz^q) + 
                 \alpha_2 y(z^qx-zx^q) + \alpha_3 z (x^qy -xy^q)\\
    &= xyz G(x,y,z),
\end{align*}
where
\[
G(x,y,z) = (\alpha_3 - \alpha_2)x^{q-1}
              +(\alpha_1 - \alpha_3)y^{q-1}
                  + (\alpha_1 - \alpha_3)z^{q-1}.
\]
The equation $G=0$ defines 
a curve described in Example~\ref{affinefilling}
by the assumption on $\alpha_1, \alpha_2$ and $\alpha_3$.

In Case~3.1,
note that
\[
A \sim 
\begin{pmatrix}
\alpha& 1&0 \\
0 &\alpha&0 \\
0 &0 &\beta
\end{pmatrix}
\sim
\begin{pmatrix}
0& 1&0 \\
0 &0&0 \\
0 &0 &\beta -\alpha
\end{pmatrix}.
\]
Put $\beta' =\beta -\alpha$, which is nonzero element of $\mathbb{F}_q$
by the assumption on $f_A(t)$.
Hence we may assume that
\begin{align*}
F_A(x,y,z) &= x(z^qx - zx^q) + \beta'z (x^qy -xy^q) \\
           & xz(z^{q-1}x - x^q + \beta'(x^{q-1} y - y^q)).
\end{align*}
Replacing $\beta'y$ by the new variable $y$, $-x$ by new $x$,
we can transform $F_A$ into
\[
-xz(x^q-xz^{q-1} + x^{q-1} y - y^q).
\]

In Case~3.2,
since
\[
A \sim 
\begin{pmatrix}
\alpha& 0&0 \\
0 &\alpha&0 \\
0 &0 &\beta
\end{pmatrix}
\sim
\begin{pmatrix}
0& 0&0 \\
0 &0&0 \\
0 &0 &\beta -\alpha
\end{pmatrix},
\]
we may assume that
\[
F_A(x,y,z) = (\beta -\alpha)z(x^qy - xy^q).
\]
So $C_A$ consists of $q+2$ $\mathbb{F}_q$-lines so that
$q+1$-lines pass through $(0,0,1)$.

In Case~4.1,
since
\[
A \sim 
\begin{pmatrix}
\alpha& 1&0 \\
0 &\alpha&1 \\
0 &0 &\alpha
\end{pmatrix}
\sim
\begin{pmatrix}
0& 1&0 \\
0 &0&1 \\
0 &0 &0
\end{pmatrix},
\]
we may assume that
\[
F_A(x,y,z) =
x(z^qx-zx^q) + y(x^qy -xy^q) \\
= xG(x,y,z).
\]
Here
$G(x,y,z) = z^qx - zx^q + x^{q-1}y^2 - y^{q+1}$.
Then the equation $-G(-z, x, y)=0$ is exactly (\ref{qplus1}).

In Case~4.2,
since
\[
A \sim 
\begin{pmatrix}
\alpha& 1&0 \\
0 &\alpha&0 \\
0 &0 &\alpha
\end{pmatrix}
\sim
\begin{pmatrix}
0& 1&0 \\
0 &0&0\\
0 &0 &0
\end{pmatrix},
\]
we may assume that
\[
F_A(x,y,z) = x(z^qx - zx^q) = x^2(z^q - zx^{q-1})
= x^2z\prod_{\lambda \in \mathbb{F}_q^{\ast}}(z- \lambda x),
\]
which shows that $C_A$
splits into $q$ reduced lines
and one double line, each of which passes through $(0,1,0)$.

In Case~4.3,
since
\[
A \sim 
\begin{pmatrix}
\alpha& 0&0 \\
0 &\alpha&0 \\
0 &0 &\alpha
\end{pmatrix},
\]
$F_A$ is the zero polynomial by Lemma~\ref{kernel}.
\qed

\begin{corollary}
For a reducible plane filling curve $C_A$ of degree $q+2$,
it contains a non-linear irreducible component if and only if
the characteristic polynomial of $A$ is the minimal polynomial.
\end{corollary}
\section{Affine-plane filling projective curve}
An affine-plane filling curve $C$ is a projective plane curve, which is not necessary irreducible, over $\mathbb{F}_q$ with the property
$C(\mathbb{F}_q)= \mathbb{A}^2(\mathbb{F}_q)$,
where $\mathbb{A}^2$ is the affine piece $\mathbb{P}^2 \setminus \{z = 0\}.$ 
For such curves, the following facts are known \cite[Theorem 3.5]{hom-kim2018}.
\begin{theorem}\label{affineplanefilling}
The degree of an affine-plane filling curve is at least $q+1$.
For a projective plane curve of degree $q+1$
in $\mathbb{P}^2$ over $\mathbb{F}_q$,
it is an affine-plane filling curve if and only if
defined by an equation of the following type:
\begin{equation}\label{equationqplusone}
(x^q-xz^{q-1}, y^q-yz^{q-1})
\left(
 \begin{array}{ccc}
   a_0 & a_1 & a_2 \\
   b_0 & b_1 & b_2
 \end{array}
\right)
\left(
\begin{array}{c}
x \\ y \\ z
\end{array}
\right) = 0,
\end{equation}
where $a_0, \dots , b_2 \in \mathbb{F}_q$
and the polynomial
\begin{equation}\label{additionalcondition}
(s,t)
\left(
 \begin{array}{cc}
   a_0 & a_1 \\
   b_0 & b_1 
 \end{array}
\right)
\left(
\begin{array}{c}
s \\ t
\end{array}
\right)
\end{equation}
in $s$ and $t$ is irreducible over $\mathbb{F}_q$.
Furthermore, an affine-plane filling curve  of degree $q+1$
is absolutely irreducible and has
a unique singular point.
\end{theorem}
Here we will add a fact supplementary to this theorem.
The maximum number of $\mathbb{F}_q$-points of a curve of degree $q+1$
without $\mathbb{F}_q$-linear components is $q^2+1$.
So the affine-plane filling curve of degree $q+1$, which is described in the above theorem, attains the second maximum.
Next theorem assures us that for not small $q$,
only the affine-plane filling curve over $\mathbb{F}_q$ of degree $q+1$
attains the second maximum for the curves of degree $q+1$.
\begin{theorem}\label{secondmaximumcurve}
Let $C$ be a plane curve of degree $q+1$ over $\mathbb{F}_q$,
which may have $\mathbb{F}_q$-linear components.
If $N_q(C)=q^2$ and $q>5$,
the missing points $\mathbb{P}^2(\mathbb{F}_q) \setminus C(\mathbb{F}_q)$
are collinear.
\end{theorem}
\proof
Let $\mathbb{P}^2(\mathbb{F}_q) \setminus C(\mathbb{F}_q)
=\{ P_0, P_1, \dots , P_{q}\}$.
Suppose that no $3$ of the $q+1$ missing points
$\{ P_0, P_1, \dots , P_{q}\}$ are collinear.
If $q$ is odd, take $(q-1)/2$ $\mathbb{F}_q$-lines
$\overline{P_{2j}P_{2j+1}}$ $(j = 1, 2, \dots , (q-1)/2)$,
where $\overline{PQ}$ is the line passing through the assigned two points
$P$ and $Q$. Also choose an $\mathbb{F}_q$-line $L_1$
such that $L_1 \ni P_1$ and $L_1 \not\ni P_0$.
Then the curve
\[
D:=L_1 \cup ( \cup_{j=1}^{\frac{q-1}{2}} \overline{P_{2j}P_{2j+1}} ) \cup C
\]
is of degree $\frac{3}{2}(q+1)$
and $D(\mathbb{F}_q) = \mathbb{P}^2(\mathbb{F}_q) \setminus \{P_0\}.$
Since $q>5$, $\frac{3}{2}(q+1)< 2q-1,$
which contradicts \cite[Corollary~2.4]{hom-kim2018}.
If $q$ is even, then
$D:=( \cup_{j=1}^{\frac{q}{2}} \overline{P_{2j-1}P_{2j}} ) \cup C$
is of degree $\frac{3}{2}q +1 < 2q-1$
and 
$D(\mathbb{F}_q) = \mathbb{P}^2(\mathbb{F}_q) \setminus \{P_0\},$
which also contradicts \cite[Corollary~2.4]{hom-kim2018}.
Therefore there are at least $3$ collinear points in the set of missing points.

Let $s:= \max
\{
|L \cap \{P_0, P_1, \dots , P_{q} \}| \mid
L \mbox{ is an $\mathbb{F}_q$-line}
\}
$
and $L_0$ an $\mathbb{F}_q$-line which attains the $s$.
Since the curve $D=L_0 \cup C$ is of degree $q+2$
and $s \geq 3$,
the condition
$D(\mathbb{F}_q) \neq \mathbb{P}^2(\mathbb{F}_q)$
implies $s=3$ by
\cite[Lemma~3.2]{hom-kim2018}.
On the other hand, from the latter part of that lemma,
the missing $q-2 (=q+1-s)$ points for $D(\mathbb{F}_q)$ should be collinear.
By the maximality of $s(=3)$, $q-2 \leq 3$, which contradicts the assumption $q>5$.
Therefore
$D(\mathbb{F}_q) = \mathbb{P}^2(\mathbb{F}_q)$,
which means the $q+1$ points $P_0, \dots , P_q$
are collinear.
\qed

\section{Fragment of affine-plane filling projective curves}
Even if the irreducible condition for (\ref{additionalcondition}) fails, 
the equation (\ref{equationqplusone}) still defines a curve $C$, which is possibly reducible,
such that $C(\mathbb{F}_q) \supset \mathbb{A}^2(\mathbb{F}_q)$.
We will investigate those curves.
\begin{lemma}\label{BtransLemma}
Let $\sigma$ be a projective transformation of $\mathbb{P}^2$
over $\mathbb{F}_q$
such that $\sigma(\{z=0\}) = (\{z=0\})$.
Denote $\sigma^{-1}\left(
\begin{array}{c}
x \\ y \\ z
\end{array}
\right) $
by
$
\left(
\begin{array}{c}
x' \\ y' \\ z'
\end{array}
\right).
$
Then there are $B \in GL(2, \mathbb{F}_q)$,
$\mathbf{b} \in (\mathbb{F}_q)^{2}$
and $\lambda \in \mathbb{F}_q^{\ast}$
such that
\begin{equation}\label{Btrans}
\left(
\begin{array}{c}
x \\ y \\ z
\end{array}
\right)
=
\left(
\begin{array}{cc}
B& \mathbf{b} \\
0& \lambda
\end{array}
\right)
\left(
\begin{array}{c}
x' \\ y' \\ z'
\end{array}
\right).
\end{equation}
If
$
F=
(x^q-xz^{q-1}, y^q-yz^{q-1})
\left(
\begin{array}{ccc}
a_0 & a_1 & a_2 \\
b_0 & b_1 & b_2
\end{array}
\right)
\left(
\begin{array}{c}
x \\ y \\ z
\end{array}
\right),
$
then $\sigma^{-1}(\{F=0 \})$
is given by the equation
\[
(x^q-xz^{q-1}, y^q-yz^{q-1})
{}^t\!B
\left(
\begin{array}{ccc}
a_0 & a_1 & a_2 \\
b_0 & b_1 & b_2
\end{array}
\right)
\left(
\begin{array}{cc}
B& \mathbf{b} \\
0& \lambda
\end{array}
\right)
\left(
\begin{array}{c}
x \\ y \\ z
\end{array}
\right)
=0.
\]
\end{lemma}
\proof
Let $B= (\beta_{i,j})$
and $\mathbf{b} = 
\left(
\begin{array}{c}
\gamma_1\\ \gamma_2
\end{array}
\right).
$
Using the relation (\ref{Btrans}),
we have
\begin{align*}
x^q-xz^{q-1} &= (\beta_{11}x' + \beta_{12}y' + \gamma_1z')^q
       - (\beta_{11}x' + \beta_{12}y' + \gamma_1z')(\lambda z')^{q-1} \\
       &= \beta_{11}(x'^q-x'z'^{q-1}) + \beta_{12}(y'^q-y'z'^{q-1}),
\end{align*}
because $\beta_{ij}^q = \beta_{ij}$
and $\lambda^{q-1}=1$.
Similarly, we have
\[
y^q-yz^{q-1} = \beta_{21}(x'^q-x'z'^{q-1}) + \beta_{22}(y'^q-y'z'^{q-1}).
\]
Summing up
$
(x^q-xz^{q-1}, y^q-yz^{q-1})
=
(x'^q-x'z'^{q-1}, y'^q-y'z'^{q-1})
{}^t\!B.
$
\qed

It is obvious that
the set
$
\left\{
\left(
\begin{array}{cc}
B& \mathbf{b} \\
0& \lambda
\end{array}
\right)
\mid
B \in GL(2, \mathbb{F}_q),\, 
\mathbf{b} \in (\mathbb{F}_q)^{2}, \, 
\lambda \in \mathbb{F}_q^{\ast}
\right\}
$
forms a subgroup of $GL(3, \mathbb{F}_q)$.
This subgroup will be denoted by $\mathcal{B}$ simply.

For
$M=
\left(
\begin{array}{ccc}
a_0&a_1&a_2\\
b_0&b_1&b_2
\end{array}
\right)
\in \mathbb{M}(2\times 3, \mathbb{F}_q),
$
the matrix consisting of first two columns
$
\left(
\begin{array}{cc}
a_0&a_1\\
b_0&b_1
\end{array}
\right)
$
is frequently denoted by $M'$
and the third column
$
\left(
\begin{array}{c}
a_2\\
b_2
\end{array}
\right)
$
by
$\mathbf{m}$.
Let
\[
G_M(x,y,z) :=
(x^q-xz^{q-1}, y^q-yz^{q-1})
\left(
 \begin{array}{ccc}
   a_0 & a_1 & a_2 \\
   b_0 & b_1 & b_2
 \end{array}
\right)
\left(
\begin{array}{c}
x \\ y \\ z
\end{array}
\right)
\]
and
\[
V_M :=\{ G_M(x,y,z)=0 \} \subset \mathbb{P}^2.
\]
From Lemma~\ref{BtransLemma},
the subgroup $\mathcal{B}\subset GL(3, \mathbb{F}_q)$ acts on
$\{ V_M \mid M \in \mathbb{M}(2\times 3, \mathbb{F}_q) \}$
as projective transformations.
For $M = (M', \mathbf{m})$ and $N = (N', \mathbf{n})$,
the symbol $M \approx N$ means that
there is a transformation
$\sigma = \left(
\begin{array}{cc}
B& \mathbf{b} \\
0& \lambda
\end{array}
\right)
\in \mathcal{B}$
such that
$N = {}^t\!B M \left(
\begin{array}{cc}
B& \mathbf{b} \\
0& \lambda
\end{array}
\right), $
that is,
$
N' = {}^t\!B M' B \mbox{\ and\ }
\mathbf{n} = {}^t\!B(M'\mathbf{b} + \lambda \mathbf{m}).
$
We also use the notation
$\sigma (M)$ as $N$ under the above situation.
Obviously the relation $\approx$ is an equivalence relation,
and $M \approx N$ if and only if
there is a projective transformation $\sigma$ over $\mathbb{F}_q$
with $\sigma(\{z=0\})=\{z=0\}$
such that $\sigma : V_M \xrightarrow{\sim}  V_N.$
In this case, we shortly say that $V_M$ is $\mathcal{B}$-equivalent to $V_N$.

We need more notation.
For $M = (M', \mathbf{m})$,
we consider the roots of the quadratic polynomial
$g_{M'}=g_{M'}(s,t) = (s, t)M'\left(
\begin{array}{c}
s\\
t
\end{array}
\right)$
in $\mathbb{P}^1$.
\begin{remark}
If $M \approx N$, then
the following statements are obviously true:
 \begin{enumerate}[(i)]
   \item $\det M' \neq 0$ if and only if $\det N' \neq 0$;
   \item $\rank M=\rank N$;
   \item $g_{M'}$ is irreducible over $\mathbb{F}_q$
   if and only if so is $g_{N'}$;
   \item $g_{M'}=0$ has two distinct roots in $\mathbb{P}^1(\mathbb{F}_q)$
   if and only if so does $g_{N'}=0$;
   \item $g_{M'}=0$ has a double root
   if and only if so does $g_{N'}=0$;
   \item $g_{M'}$ is the zero polynomial if and only if so is $g_{N'}$.
 \end{enumerate}
\end{remark}
\begin{remark}
If $\det M' = 0$, then
$g_{M'}$ is either reducible or the zero polynomial.
The converse is not always true.
\end{remark}
\proof
Let $M'= \left(
 \begin{array}{cc}
   a_0 & a_1 \\
   b_0 & b_1 
 \end{array}
\right).$
If $a_0=b_0=0$,
then $g_{M'} =t(a_1s +b_1t)$.
If $(a_0, b_0) \neq (0,0)$, there is an element $\lambda \in \mathbb{F}_q$
such that $(a_1, b_1) = \lambda (a_0, b_0)$.
Hence
$
g_{M'} = (a_0 s + b_0t)(s + \lambda t).
$
A matrix $M'$
with $b_0 =-a_1$ and $b_1=0$
gives an example for which the converse does not hold.
\qed

\begin{proposition}\label{listofM}
Let $M=(M', \mathbf{m})$ be a nonzero $2\times 3$ matrix over $\mathbb{F}_q$
such that $g_{M'}$ is reducible or the zero polynomial.
\begin{enumerate}[{\rm (I)}]
 \item Suppose the equation $g_{M'}=0$ has two distinct roots.
    \begin{enumerate}[{\rm (\mbox{I-}1)}]
      \item If $\det M' \neq 0$, then
       $
       M \approx 
       \left(
 \begin{array}{ccc}
   0 & a_1 & 0 \\
   b_0 &0 & 0
 \end{array}
\right)
       $
       with $a_1b_0 \neq 0$ and $a_1 +b_0 \neq 0$.
      \item If $\det M' = 0$ and $\rank M=2$, then
      $
       M \approx 
       \left(
 \begin{array}{ccc}
   0 & a_1 & 0 \\
   0 &0 & b_2
 \end{array}
\right)
       $
       with  $a_1b_2 \neq 0$.
       \item If $\det M' = 0$ and $\rank M=1$, then
      $
       M \approx 
       \left(
 \begin{array}{ccc}
   0 & a_1 & 0 \\
   0 &0 & 0
 \end{array}
\right)
       $
       with  $a_1\neq 0$.
    \end{enumerate}
  \item Suppose the equation $g_{M'}=0$ has a double root.
      \begin{enumerate}[{\rm (\mbox{II-}1)}]
      \item
       If $\det M' \neq 0$, then
       $
       M \approx 
       \left(
 \begin{array}{ccc}
   a_0 & a_1 & 0 \\
   -a_1&0 & 0
 \end{array}
\right)
       $
       with $a_0a_1 \neq 0$.
      \item
      If $\det M' = 0$ and $\rank M=2$,
      then
      $
       M \approx 
       \left(
 \begin{array}{ccc}
   a_0 & 0 & 0 \\
   0 &0 & b_2
 \end{array}
\right)
       $
       with $a_0b_2 \neq 0$.
      \item
      If $\det M'=0$ and $\rank M =1$,
      $
       M \approx 
       \left(
 \begin{array}{ccc}
   a_0 & 0 & 0 \\
   0 &0 & 0
 \end{array}
\right)
       $
       with $a_0 \neq 0$.
    \end{enumerate}
   \item Suppose $g_{M'}$ is the zero polynomial.
      \begin{enumerate}[{\rm (\mbox{III-}1)}]
      \item
      If $\det M' \neq 0$, then
      $
       M \approx 
       \left(
 \begin{array}{ccc}
   0 & a_1 & 0 \\
   -a_1 &0 & 0
 \end{array}
\right)
       $
       with $a_1 \neq 0$.
      \item
      There is no member $M$ so that $g_{M'}$ is the zero polynomial,
      $\det M'=0$ and $\rank M =2$.
      \item
      If $\det M' = 0$ and $rank~M =1$, then
      $
       M \approx 
       \left(
 \begin{array}{ccc}
   0 & 0 & 1 \\
   0 &0 & 0
 \end{array}
\right)
       $.
    \end{enumerate}
\end{enumerate}
\end{proposition}

This proposition is just summary of the following series of lemmas.
\begin{lemma}\label{lemma-1}
If $\det M' \neq 0$, then $M \approx (M', \mathbf{0})$.
\end{lemma}
\proof
For $M=(M', \mathbf{m})$,
there is a vector $\mathbf{b} \in \mathbb{F}_q^2$
such that $M' \mathbf{b} = -\mathbf{m}$.
Choose $\sigma \in \mathcal{B}$
as $\sigma = \left(
 \begin{array}{cc}
   1_2 & \mathbf{b} \\
   0 &1
 \end{array}
\right).
$
Then $\sigma (M)= (M', \mathbf{0}).$
\qed
\begin{lemma}\label{lemma-2}
If $\det M' \neq 0$
and the equation $g_{M'}=0$ has two distinct roots, then
       $
       M \approx 
       \left(
 \begin{array}{ccc}
   0 & a_1 & 0 \\
   b_0 &0 & 0
 \end{array}
\right)
       $
       with $a_1b_0 \neq 0$ and $a_1 +b_0 \neq 0$.
\end{lemma}
\proof
Let $(\alpha_1, \alpha_2)$ and  $(\beta_1, \beta_2)$
be two distinct roots of $g_{M'}(s,t)=0$.
Then 
$
\det \left(
 \begin{array}{cc}
   \alpha_1 & \beta_1 \\
   \alpha_2 & \beta_2
 \end{array}
\right) \neq 0
$
and
$
{}^{t}\!\left(
 \begin{array}{cc}
   \alpha_1 & \beta_1 \\
   \alpha_2 & \beta_2
 \end{array}
\right)M'
\left(
 \begin{array}{cc}
   \alpha_1 & \beta_1 \\
   \alpha_2 & \beta_2
 \end{array}
\right)
=
\left(
 \begin{array}{cc}
   0 & a_1 \\
   b_0 & 0
 \end{array}
\right)
$
for some $a_1, b_0 \in \mathbb{F}_q$.
Hence
$M
\approx
\left(
 \begin{array}{ccc}
   0 & a_1 &0\\
   b_0 & 0 &0
 \end{array}
\right)
$
by Lemma~\ref{lemma-1}.
The two properties on $M'$ are hereditary under the relation $\approx$.
\qed
\begin{lemma}\label{lemma-3}
If $\det M' = 0$
and the equation $g_{M'}=0$ has two distinct roots, then
       $
       M \approx 
       \left(
 \begin{array}{ccc}
   0 & a_1 & 0 \\
   0 &0 & b_2
 \end{array}
\right)
       $
       for some $a_1 , \,  b_2 \in \mathbb{F}_q$ with $a_1 \neq 0$.
\end{lemma}
\proof
Let $(\alpha_1, \alpha_2)$ and  $(\beta_1, \beta_2)$
be two distinct roots,
and
$
\sigma =
\left(
 \begin{array}{ccc}
   \alpha_1 & \beta_1 &0\\
   \alpha_2 & \beta_2 &0\\
           0&    0    &1
 \end{array}
\right)
\in \mathcal{B}.
$
Then
$
\sigma (M) =
\left(
 \begin{array}{ccc}
   0 & a_1' & a_2' \\
   b_0' &0 & b_2'
 \end{array}
\right)
$
for some $a_1' , \dots ,  b_2' \in \mathbb{F}_q$
and one of $\{ a_1' , b_0'\}$ is $0$ and other is nonzero,
because $\det \sigma (M)'=0$
and $g_{\sigma (M)'}=0$ has two distinct roots.
Let $\tau
   = \left(
 \begin{array}{ccc}
   0 & 1 &0\\
   1 & 0 &0\\
   0&  0 &1
 \end{array}
\right)
\in \mathcal{B}.
$
Then
$
\tau(\sigma ( M))
= \left(
 \begin{array}{ccc}
   0 &b_0' & a_2' \\
   a_1'  &0 & b_2'
 \end{array}
\right).
$
Therefore $
M \approx 
       \left(
 \begin{array}{ccc}
   0 & a_1 & a_2 \\
   0 &0 & b_2
 \end{array}
\right)
       $
       for some $a_1 , a_2,  b_2 \in \mathbb{F}_q$ with $a_1 \neq 0$.
Finally,
let
$
\rho =
\left(
 \begin{array}{ccc}
   1 & 0 &0\\
   0& 1 & -\frac{a_2}{a_1}\\
   0&  0 &1
 \end{array}
\right)
\in \mathcal{B}.
$
Then
$
\rho(M) \approx 
\left(
 \begin{array}{ccc}
   0 & a_1 & 0 \\
   0 &0 & b_2
 \end{array}
\right).
$
\qed
\begin{lemma}\label{lemma-4}
Suppose the equation $g_{M'}=0$ has a double root.
\begin{enumerate}[{\rm (1)}]
 \item If $\det M' \neq 0$, then
  $
       M \approx 
       \left(
 \begin{array}{ccc}
   a_0 & a_1 & 0 \\
   -a_1&0 & 0
 \end{array}
\right)
       $
     with $a_0a_1 \neq 0$.
 \item If $\det M' = 0$, then
  $
       M \approx 
       \left(
 \begin{array}{ccc}
   a_0 & 0 & 0 \\
   0&0 & b_2
 \end{array}
\right)
       $
     with $a_0 \neq 0$.
\end{enumerate}
\end{lemma}
\proof
Let $(\alpha_1, \alpha_2)$ be the root of $g_{M'}=0$.
Choose a vector $(\beta_1, \beta_2)$
which is not a multiple of the first choice.
Let
$B=
\left(
 \begin{array}{cc}
  \beta_1 & \alpha_1 \\
  \beta_2 &  \alpha_2 
 \end{array}
\right)
$
and
$
\sigma =
\left(
 \begin{array}{cc}
  B& \mathbf{0} \\
  {}^{t}\mathbf{0} &  1 
 \end{array}
\right).
$
Then
$
\sigma (M) = 
\left(
 \begin{array}{ccc}
   a_0 & a_1 & a_2 \\
   b_0&0 & b_2
 \end{array}
\right).
$
Since $g_{\sigma(M)'} =0$
has a double root,
$a_1+b_0 =0$ and $a_0 \neq 0$.

(1) If $\det M' \neq 0$, then $\det \sigma(M)' \neq 0$ either.
Hence $a_1 \neq 0$ and
$M \approx \sigma(M) \approx
\left(
 \begin{array}{ccc}
   a_0 & a_1 & 0 \\
   -a_1&0 & 0
 \end{array}
\right)
$
by Lemma~\ref{lemma-1}.

(2)  If $\det M' = 0$, then $a_1=b_0=0$.
Furthermore,
let
$
\tau
= \left(
 \begin{array}{ccc}
   1 & 0 &-\frac{a_2}{a_0}\\
   0& 1 & 0 \\
   0&  0 &1
 \end{array}
\right).
$
Then
$
\tau (\sigma (M)) \approx
\left(
 \begin{array}{ccc}
   a_0 & 0 & 0 \\
   0&0 & b_2
 \end{array}
\right).
$
\qed
\begin{lemma}\label{lemma-5}
Suppose $g_{M'}$ is the zero polynomial.
   \begin{enumerate}[{\rm (1)}]
    \item If $\det M' \neq 0$, then
     $
     M \approx 
       \left(
 \begin{array}{ccc}
   0 & a_1 & 0 \\
   -a_1&0 & 0
 \end{array}
\right)
     $
     with $a_1 \neq 0$.
    \item If $\det M' = 0$, then
     $
     M \approx 
       \left(
 \begin{array}{ccc}
   0 & 0 & 1 \\
   0&0 & 0
 \end{array}
\right).
     $
   \end{enumerate}
\end{lemma}
\proof
Let
$
M' = 
\left(
 \begin{array}{cc}
   a_0 & a_1 \\
   b_0& b_1
 \end{array}
\right).
$
Then $g_{M'}$ is the zero polynomial
if and only if $a_0=b_1=0$ and $b_0=-a_1$.

(1) If $\det M' \neq 0$, then $a_1 \neq 0$ and
     $
     M \approx 
       \left(
 \begin{array}{ccc}
   0 & a_1 & 0 \\
   -a_1&0 & 0
 \end{array}
\right)
     $
     by Lemma~\ref{lemma-1}.
     
 (2) If $\det M'=0$, then
 $
 M =
       \left(
 \begin{array}{ccc}
   0 & 0 & a_2 \\
   0 &0 & b_2
 \end{array}
\right).
     $
Since $M$ is a nonzero matrix,
$(a_2, b_2) \neq (0,0)$.
Hence we can find a mtrix $B \in GL(2, \mathbb{F}_q)$
such that
${}^t\!B\left(
 \begin{array}{c}
    a_2 \\
    b_2
 \end{array}
\right)
=
\left(
 \begin{array}{c}
    1\\
    0
 \end{array}
\right).
$
Put
$\sigma =
\left(
 \begin{array}{cc}
  B& \mathbf{0} \\
  {}^{t}\mathbf{0} &  1 
 \end{array}
\right).
$
Then
$\sigma (M)
= \left(
 \begin{array}{ccc}
   0 & 0 & 1 \\
   0&0 & 0
 \end{array}
\right).
$
\qed

\gyokan

Finally, We give a description of $V_M$ according to the list in
Proposition~\ref{listofM}.

Note that if $M \approx N$, then
$
|V_M(\mathbb{F}_q) \cap \{z=0\}| = |V_N(\mathbb{F}_q)  \cap \{z=0\}|,
$
which is called the number of $\mathbb{F}_q$-points at infinity of $V_M$.

\begin{theorem}
Let $M=(M', \mathbf{m}) \in \mathbb{M}(2\times 3, \mathbb{F}_q)$
such that $g_{M'}$ is reducible or the zero polynomial.
\begin{enumerate}[{\rm (I)}]
 \item If $g_{M'}=0$ has two distinct roots,
 then the number of $\mathbb{F}_q$-points at infinity of $V_M$ is two.
    \begin{enumerate}[{\rm (\mbox{I-}1)}]
     \item In addition, 
     if $\det M' \neq 0$, then
     $V_M$ is $\mathcal{B}$-equivalent to the curve
     \[
     xy(a_1 x^{q-1}+ b_0y^{q-1} - (a_1+b_0)z^{q-1}) =0
     \]
     with $a_1b_0(a_1+b_0) \neq 0$.
     The nonlinear component is the curve described in Example~\ref{degreeqminusone}.
     \item In addition, if $\det M' = 0$ and $\rank M=2$, then
      $V_M$ is $\mathcal{B}$-equivalent to the curve
     \[
     y(x^q-xz^{q-1} + y^{q-1}z - z^q)=0.
     \]
     The nonlinear component is the curve described in Example~\ref{degreeq}.
     \item In addition, if $\det M' = 0$ and $\rank M=1$, then
$V_M$ is $\mathcal{B}$-equivalent to the curve
     \[
     y(x^q - xz^{q-1})=0,
     \]
     which consists of $q$ lines passing through $(0,1,0)$
     and another line not passing through the point.
     \end{enumerate}
\item If $g_{M'}=0$ has a double root,
 then the number of $\mathbb{F}_q$-points at infinity of $V_M$ is one.
    \begin{enumerate}[{\rm (\mbox{II-}1)}]
     \item In addition, if $\det M' \neq 0$, 
       then $V_M$ is $\mathcal{B}$-equivalent to the curve
     \[
     x(x^q -xz^{q-1} + yx^{q-1} -y^q)=0.
     \]
     The nonlinear component is the curve described in Example~\ref{degreeq}.
     \item In addition, if $\det M' = 0$ and $\rank M=2$, then
      $V_M$ is $\mathcal{B}$-equivalent to the curve
     \[
     x^{q+1}-x^2z^{q-1}+y^qz-yz^q=0,
     \]
     which is nonsingular and described in Example~\ref{degreeqplusone}.
     \item In addition, if $\det M' = 0$ and $\rank M=1$, then
      $V_M$ is $\mathcal{B}$-equivalent to the curve
     \[
     x^2(x^{q-1} -z^{q-1})=0,
     \]
     which consists of $q-1$ simple lines and one double line.
     Those lines pass through a common point.
    \end{enumerate}
\item If $g_{M'}$ is the zero polynomial, then $V_M(\mathbb{F}_q) = \mathbb{P}^2(\mathbb{F}_q)$ if $V_M$ exists.
     \begin{enumerate}[{\rm (\mbox{III-}1)}]
     \item In addition, if $\det M' \neq 0$, 
       then $V_M$ is $\mathcal{B}$-equivalent to the curve
     \[
     x^{q}y- xy^{q}=0,
     \]
     which is the union of the $q+1$ $\mathbb{F}_q$-lines passing through $(0,0,1)$.
     \item There is no curve for $M$ with properties $g_{M'}$ is the zero polynomial, $\det M' = 0$ and $\rank M=2$.
     \item  In addition, if $\det M' = 0$ and $\rank M=1$, then
     $V_M$ is $\mathcal{B}$-equivalent to the curve
     \[
     x^{q}z- xz^{q}=0,
     \]
     which is the union of the $q+1$ $\mathbb{F}_q$-lines passing through $(0,1,0)$.
     \end{enumerate}
\end{enumerate}
\end{theorem}
\proof
Since $V_M \cap \{z=0 \}$ is defined by
$G_M(x,y,0)=
(x^q, y^q)M'\left(
 \begin{array}{c}
   x\\
   y
 \end{array}
\right)
=0$,
$V_M(\mathbb{F}_q) \cap \{z=0\}
= \{ (\alpha, \beta, 0) \in \mathbb{P}^2(\mathbb{F}_q)
\mid g_{M'}(\alpha, \beta) =0 \}$.

Straightforward computations from ``canonical" forms  listed in Proposition~\ref{listofM} induce the equation in each case.
\qed
\begin{corollary}
For $M=(M', \mathbf{m}) \in \mathbb{M}(2\times 3, \mathbb{F}_q)$,
$V_M$ has a nonlinear irreducible component if and only if
$\rank M =2$ and
$g_{M'}$ is a nonzero polynomial.
\end{corollary}



\begin{thebibliography}{00}
\bibitem{dur2018}
G. Duran Cunha,
{\em Curves containing all points of a finite projective Galois plane},
J. Pure Appl. Algebra 222 (2018) 2964--2974.
\bibitem{gan}
F. R. Gantmacher,
The theory of matrices. Vol. 1
(Translated from the Russian by K. A. Hirsch. Reprint of the 1959 translation.)
AMS Chelsea Publishing, Providence, 1998. 
\bibitem{hom-kim2010b}
M. Homma and S. J. Kim,
{\em Sziklai's conjecture on the number of points of
a plane curve over a finite field  {\rm III}},
Finite Fields Appl. 16 (2010) 315--319.
\bibitem{hom-kim2011}
M. Homma and S. J. Kim,
{\em Toward determination of optimal plane curves
with a fixed degree over a finite field}, 
Finite Fields Appl. 17 (2011) 240--253.
\bibitem{hom-kim2012}
M. Homma and S. J. Kim,
{\em The uniqueness of a plane curve of degree $q$
attaining Sziklai's bound over $\mathbb{F}_q$}, 
Finite Fields Appl. 18 (2012) 567--580.
\bibitem{hom-kim2013}
M. Homma and S. J. Kim,
{\em Nonsingular plane filling curves of minimum degree
over a finite field and their automorphism groups:
Supplements to a work of Tallini},
Linear Algebra Appl. 438 (2013) 969--985.
\bibitem{hom-kim2018}
M. Homma and S. J. Kim,
{\em The second largest number of points on plane
curves over finite fields},
Finite Fields Appl. 49 (2018) 80--93.
\bibitem{szi2008}
P. Sziklai,
{\em A bound on the number of points of a plane curve}, 
Finite Fields Appl. 14 (2008) 41--43.
\bibitem{tal1}
G. Tallini,
{\em Le ipersuperficie irriducibili d'ordine minimo che
invadono uno spazio di Galois},
Atti Accad. Naz. Lincei Rend. Cl. Sci. Fis. Mat. Nat. (8) 30 (1961) 706--712.
\bibitem{tal2}
G. Tallini,
{\em Sulle ipersuperficie irriducibili d'ordine minimo che
contengono tutti i punti di uno spazio di Galois $S_{r,q}$},
Rend. Mat. e Appl. (5) 20 (1961) 431--479.
\end{thebibliography}
\end{document}